\documentclass[12pt,dvips]{amsart}
\usepackage{euler, amsfonts, amssymb, latexsym, epsfig}
\renewcommand{\familydefault}{ppl}

\newcommand\B{\reals^3_{\sum=0}}
\newcommand\BZ{\integers^3_{\sum=0}}

\newcommand\Tensor{{\bigotimes}}
\newcommand\calO{{\mathcal O}}
\newcommand\calM{{\mathcal M}_{\lambda\mu\nu}}
\newtheorem{Theorem}{Theorem} 
\newtheorem{conjecture}{Conjecture} 
 
\newtheorem{Lemma}{Lemma}

\newtheorem*{Theorem*}{Theorem}

\newcommand\reals{{\mathbb R}}
\newcommand\complexes{{\mathbb C}}
\newcommand\integers{{\mathbb Z}}

\newcommand\Ldots{\ldots\hskip-.03cm}
\newcommand\eps{\varepsilon}
\theoremstyle{plain}

\newcommand\dfn{\bf} % maybe should be \em

\newcommand\HONEY{{\tt HONEY}}
\newcommand\BDRY{{\tt BDRY}}

\begin{document}
\pagestyle{plain}

\title{Honeycombs and sums of Hermitian matrices}
\author{Allen Knutson}
\thanks{AK's research was partially conducted for the
Clay Mathematics Institute.}
\email{allenk@math.berkeley.edu}
\address{Mathematics Department\\ UC Berkeley\\ Berkeley, California}
\author{Terence Tao}
\thanks{TT is supported by the Clay Mathematics Institute and by
  grants from the Sloan and Packard foundations.}
\email{tao@math.ucla.edu}
\address{Mathematics Department\\ UCLA\\ Los Angeles, California}
\date{\today}

\begin{abstract}
  Horn's conjecture \cite{Horn}, which given the spectra of two
  Hermitian matrices describes the possible spectra of the sum,
  was recently settled in the affirmative.
  We discuss one of the many steps in this, which
  required us to introduce a combinatorial gadget called 
  a {\em honeycomb}; the question is then reformulable as about the
  existence of honeycombs with certain boundary conditions.  Another
  important tool is the connection to the representation theory of
  the group $U(n)$, by ``classical vs.  quantum'' analogies.
\end{abstract}

\maketitle

%\tableofcontents\label{fig:}

If $H$ is an $n \times n$ Hermitian matrix, then the spectrum $\lambda = (\lambda_1 \geq \ldots \geq \lambda_n)$ can be written as a weakly decreasing sequence of $n$ real numbers.  Conversely, for every spectrum $\lambda$ we can form the set $\calO_\lambda$ of Hermitian matrices with spectrum $\lambda$; this set is known as a \emph{co-adjoint orbit} of $U(n)$.

If $\lambda$, $\mu$, $\nu$ are three spectra, we define the relation
\begin{equation}\label{class-eq}
\lambda \boxplus \mu \sim_c \nu
\end{equation}
if there exist Hermitian matrices $H_\lambda \in \calO_\lambda$, $H_\mu \in \calO_\mu$, $H_\nu \in \calO_\nu$ such that $H_\lambda + H_\mu = H_\nu$.  (The ``c'' in $\sim_c$ stands for ``classical''; we will define a quantum analogue $\sim_q$ later on.)  For instance, one can easily verify that 
$$(3) \boxplus (4) \sim_c (7), \quad (3,0) \boxplus (4,0) \sim_c (7,0), \quad
(3,0) \boxplus (4,0) \sim_c (4,3), \quad (2,0) \boxplus (2,0) \sim_c (3,1)$$
but that
$$(3) \boxplus (4) \not \sim_c (5), \quad (3,0) \boxplus (4,0) \not \sim_c (8, -1)$$

In 1912, Hermann Weyl \cite{Weyl} posed the problem of determining the set of triples $\lambda$, $\mu$, $\nu$ for which \eqref{class-eq} held.  Or more informally: given the eigenvalues of two Hermitian matrices $H_\lambda$ and $H_\mu$, what are all the possible eigenvalues of the sum $H_\lambda + H_\mu$?
The purpose of this article is to describe the successful resolution to this problem, based on recent breakthroughs \cite{Klyachko,HR,KT,KTW}.

It is fairly easy to obtain necessary conditions in order for \eqref{class-eq} to hold.  For instance, from the simple observation that the trace of $H_\lambda + H_\mu$ must equal the sum of the traces of $H_\lambda$ and $H_\mu$, we obtain the condition
\begin{equation}\label{trace}
\nu_1 + \ldots + \nu_n = \lambda_1 + \ldots + \lambda_n + \mu_1 + \ldots + \mu_n.
\end{equation}
Another immediate constraint is that 
\begin{equation}\label{111}
\nu_1 \leq \lambda_1 + \mu_1,
\end{equation}
since the largest eigenvalue of $H_\lambda + H_\mu$ is at most the sum of $H_\lambda$ and $H_\mu$'s individual largest
eigenvalues. (Exercise for the reader: show equality occurs exactly when
the same vector is a principal eigenvector for both matrices.)
Weyl found a number of similar necessary conditions, such as
the statement $\nu_{i+j+1} \leq \lambda_{i+1} + \mu_{j+1}$ whenever $0 \leq i, j, i+j < n$, and many more necessary conditions were found by later authors.  All of these conditions took the form of homogeneous linear inequalities (e.g. $\nu_1 + \nu_2 \leq \lambda_1 + \lambda_2 + \mu_1 + \mu_2$.)  This phenomenon was then explained in the work of Alfred Horn \cite{Horn}.  Based on extensive computations for small $n$, Horn then described a recursive algorithm for generating a list of such linear inequalities, and conjectured that this list, together with the trace identity \eqref{trace}, gave a necessary and sufficient set of criteria for \eqref{class-eq}. His conjecture turned out
to be correct, though it waited 36 years for resolution.

\begin{figure}[htbp]
  \begin{center}
    \leavevmode
    \epsfig{file=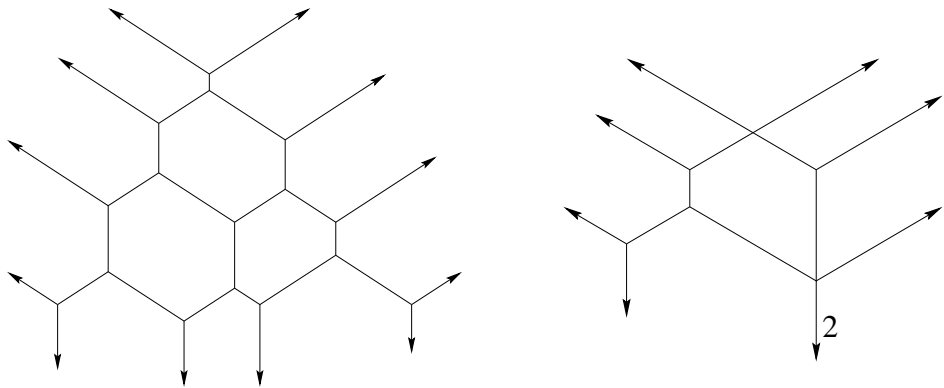,height=2in}
    \caption{Two honeycombs. The left one is more typical in having only
      Y vertices, as will be explained in theorem \ref{thm:picturesdontlie}.
      All edges are multiplicity 1, except for the edge labeled 2 
      in the right-hand honeycomb.}
    \label{fig:honeyex}
  \end{center}
\end{figure}

We approached this problem by first observing that the relation \eqref{class-eq} could be converted into a statement about {\em honeycombs}, which we introduced (for this purpose) in \cite{KT}. These are a family of planar arrangements of edges labeled with multiplicities (some examples are in figure \ref{fig:honeyex}.) We give the precise definition in the next section.

The relevance of honeycombs to sums of Hermitian matrices is the following theorem, which we explain in more detail later:

\begin{Theorem}\label{thm:cexist}
The relationship \eqref{class-eq} holds if and only if
there exists a honeycomb with boundary values $(\lambda,\mu,-\nu)$.  
\end{Theorem}

We now briefly summarize the sequence of events leading to the full proof of Horn's conjecture.  It has been known for some time that the ``classical'' problem \eqref{class-eq} of summing Hermitian matrices is closely related to the ``quantum'' problem 
\begin{equation}\label{quantum-eq}
\lambda \boxplus \mu \sim_q \nu
\end{equation} 
of tensoring $U(n)$ representations, which we shall define later.  For instance, it was well-known that \eqref{quantum-eq} implies \eqref{class-eq}.  The converse implication (for integral $\lambda, \mu, \nu$) was known as the \emph{saturation conjecture}.
Using some formidable algebraic and geometric machinery, Klyachko \cite{Klyachko} was able to demonstrate a further non-trivial recursive relationship between the classical and quantum problems, and noted that this, combined with the saturation conjecture, would imply Horn's conjecture.  In \cite{KT}, Theorem \ref{thm:cexist} (and the quantum analogue of this theorem) were used to convert the saturation conjecture into a statement about honeycombs, and then proved this statement by combinatorial methods, thus proving the saturation and Horn conjectures. (Recently it has been shown \cite{KTW} that one can derive Horn's conjecture directly from the saturation conjecture by purely combinatorial techniques, bypassing the machinery of \cite{Klyachko}. Also, a very different proof of saturation, based on the representation theory of quivers, has since been given in \cite{DW}.  Finally, a short rendition of \cite{KT} can be found in \cite{Bu}.)

We shall give a rather ahistorical (and pro-honeycomb) survey of this circle of ideas, starting with honeycombs (which were actually the last piece of the puzzle to be discovered), then discussing the connections between the classical and quantum problems, followed by a sketch of the honeycomb-based proof of the saturation conjecture.  Then we state Horn's conjecture, and sketch the honeycomb-based proof of this from saturation.

There are many other closely related and interesting mathematical questions  that we will not address, and we shall refer to the reader to the excellent survey article \cite{FultonBAMS}. 

\section{Honeycombs}\label{honey-sec}

We shall shortly give the formal definition of a honeycomb as introduced in \cite{KT}, but first we shall rephrase Weyl's problem in a more symmetric form.  We say that the relation
\begin{equation}\label{lmn}
\lambda \boxplus \mu \boxplus \nu \sim_c 0
\end{equation}
holds if there exist matrices $H_\lambda \in \calO_\lambda$, $H_\mu \in \calO_\mu$, $H_\nu \in \calO_\nu$ such that $H_\lambda + H_\mu + H_\nu = 0$.  Clearly we have
$$ \lambda \boxplus \mu \sim_c \nu \iff \lambda \boxplus \mu \boxplus (-\nu) \sim_c 0$$
where $-\nu := (-\nu_n, \ldots, -\nu_1)$ is the negation of $\nu$.  Thus to solve Weyl's problem it suffices to determine the set of triples $\lambda, \mu, \nu$ which obey \eqref{lmn}.  This formulation has the advantage of $S_3$ symmetry in $(\lambda, \mu, \nu)$, as opposed to mere $S_2$ symmetry in $(\lambda, \mu)$.

The trace identity \eqref{trace} in this new formulation becomes
\begin{equation}\label{trace-sym}
\lambda_1 + \ldots + \lambda_n + \mu_1 + \ldots + \mu_n + \nu_1 + \ldots + \nu_n = 0,
\end{equation}
while \eqref{111} becomes
\begin{equation}\label{11n}
\lambda_1 + \mu_1 + \nu_n \geq 0.
\end{equation}
Based on these relations, it is natural to introduce the plane
$$ \B := \{(x,y,z) \in \reals^3 : x+y+z=0\}.$$
We shall always depict this plane with the six 
``cardinal directions'' $(0,1,-1)$, $(-1, 1, 0)$, $(-1,0,1)$, $(0, -1, 1)$ and $(1,-1,0)$, and $(1, 0, -1)$ drawn Northwest, North, Northeast, Southeast, South, and Southwest respectively.  (Of course, ``Northwest'' makes a 60 degree angle with North rather than a 45 degree angle, and similarly for the other diagonal cardinal directions.) 

A \emph{diagram} will be defined as a configuration of (possibly half-infinite) line segments in $\B$, with each edge parallel to one of the cardinal directions (North-South, Northeast-Southwest, Northwest-Southeast), and labeled with a positive integer which we refer to as the ``multiplicity'' or ``tension''.  To every diagram we can associate a measure on $\B$, defined as the sum of Lebesgue measure on each line segment, weighted by the multiplicity.  We say that two diagrams $d$, $d'$ are \emph{equivalent} if their associated measures are equal.

If $h$ is a diagram and $v$ is a point in $\B$, we say that $v$ is a {\dfn zero-tension point} of $h$ if, in a sufficiently small neighbourhood of $v$, $h$ is equivalent to a union of rays emanating from $v$, and the sum of the co-ordinate vectors of these rays, multiplied by their tensions, equals zero.

\begin{figure}[htbp]
  \begin{center}
    \leavevmode
    \epsfig{file=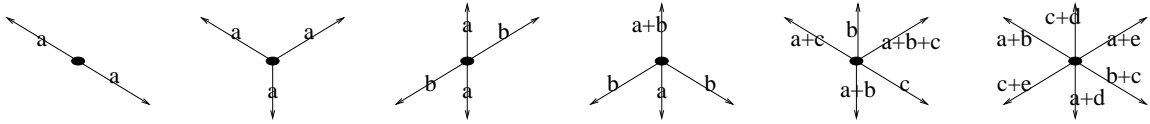,width=6in}
    \caption{Zero-tension points, with the rays labeled by their 
      multiplicities. All but the first type are called `vertices'.}
    \label{fig:ztensvs}
  \end{center}
\end{figure}

The two possibilities that will interest us most are a point on a line
segment, in which case the zero-tension condition says that the two
rays must have the same multiplicity, and a point at the center of a Y with
again three equal-multiplicity rays. There are several more
complicated possibilities, as shown in figure \ref{fig:ztensvs}.

Define a {\dfn honeycomb} $h$ as a diagram (or more precisely, an equivalence class of diagrams) such that 
\begin{enumerate}
\item every point in $\B$ is a zero-tension point
\item there are only finitely many ``vertices,'' i.e. points with more than
  two rays emanating
\item the semiinfinite lines go only in the
  Northeast, Northwest, and South directions (i.e. no Southeast, Southwest, or North rays)
\end{enumerate}

The lines mentioned in part 3 are called the {\dfn boundary edges} of
the honeycomb. Two examples of honeycombs appear in figure \ref{fig:honeyex}.

It is a pleasant exercise to show that the number of boundary edges
(with multiplicity) pointing in one cardinal direction is the same 
as the number in each of the other two directions.  (This is basically because the net tension of the honeycomb must be zero.)  We will call a honeycomb
with $n$ boundary edges in each direction an {\dfn $n$-honeycomb} and
denote the space of such by $\HONEY_n$.

\begin{figure}[htbp]
  \begin{center}
    \input{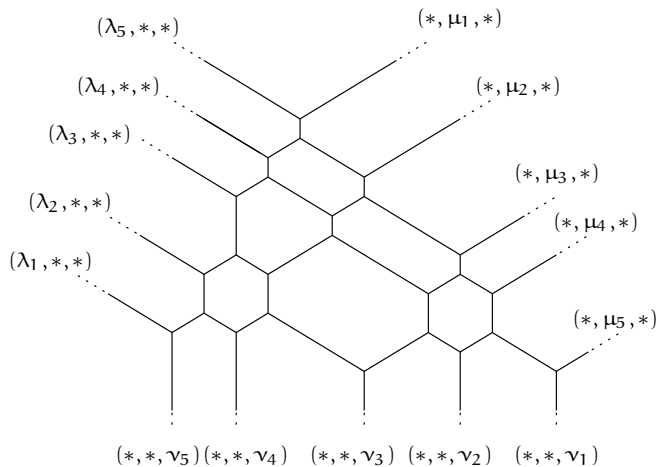}
    \caption{The constant coordinates on the boundary edges
      of a $\tau_5$-honeycomb. (The stars are the nonconstant coordinates.)}
    \label{fig:honeylabel}
  \end{center}
\end{figure}

Since every edge in a honeycomb is parallel to a cardinal direction,
each of which has one of its three coordinates being zero, every
honeycomb edge has a {\dfn constant coordinate} (common to every point
along that edge.) In particular we can read off the constant coordinates 
of boundary edges and call them 
$$(\lambda_1,\Ldots,\lambda_n,\mu_1,\Ldots,\mu_n,\nu_1,\Ldots,\nu_n) = (\lambda, \mu, \nu)$$
as in figure \ref{fig:honeylabel}. 

We can now phrase Theorem \ref{thm:cexist} in this symmetrized setting:

\begin{Theorem}\label{thm:cexist-sym}
The relationship $\lambda \boxplus \mu \boxplus \nu \sim_c 0$ holds if and only if
there exists a honeycomb with boundary values $(\lambda,\mu,\nu)$.  
\end{Theorem}

Interestingly, almost all the proofs we know of this Theorem proceed by first proving a quantized version of this Theorem, which we define in the next section. We shall therefore not discuss the proof of this Theorem here, and content ourselves instead with producing evidence which strongly suggests that the Theorem is plausible.

\begin{figure}[htbp]
  \begin{center}
\epsfig{file=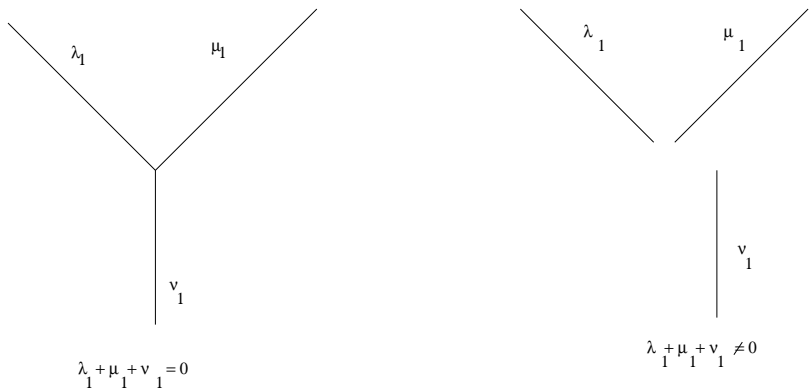,height=2in}
\caption{A 1-honeycomb can be formed if and only if the boundary values sum to zero.  The edges are labeled by their constant co-ordinates.}
    \label{fig:1dim}
  \end{center}
\end{figure}

We first consider the $n=1$ case.  In this case $\lambda = (\lambda_1)$, $\mu = (\mu_1)$, $\nu = (\nu_1)$, and it is clear that \eqref{lmn} holds if and only if $\lambda_1 + \mu_1 + \nu_1 = 0$.  (In other words, the trace condition is already necessary and sufficient.)  On the honeycomb side, this claim can be easily seen, if one accepts the fact (which is actually a little tricky to prove) that $1$-honeycombs must have the shape of a ``Y''.  See Figure \ref{fig:1dim}.

More generally, it is a pleasant exercise to show that the boundary values of any $n$-honeycomb must satisfy \eqref{trace-sym}, basically because the three co-ordinates around every vertex sum to zero (by virtue of lying in $\B$).

\begin{figure}[htbp]
  \begin{center}
\leavevmode
\epsfig{file=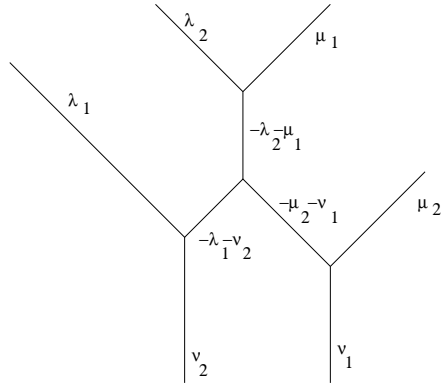,height=2in}
\caption{A 2-honeycomb is uniquely determined by its boundary values. The boundary values must satisfy \eqref{trace-sym}, and the edge lengths must be non-negative.  The edges are labeled by their constant co-ordinates.  The edge lengths can be computed by subtracting the constant co-ordinates of two parallel adjacent edges; for instance, the lower left edge has length $\lambda_1 - (-\mu_2 - \nu_1)$.}
    \label{fig:2dim}
  \end{center}
\end{figure}

Now consider the $n=2$ case, so that $\lambda = (\lambda_1, \lambda_2)$, $\mu = (\mu_1, \mu_2)$, $\nu = (\nu_1, \nu_2)$.  In this case there can be at most one 2-honeycomb with the specified boundary values (Figure \ref{fig:2dim}).  The lengths of the three line segments in the honeycomb can be computed as $\lambda_2 + \mu_1 + \nu_1$, $\lambda_1 + \mu_2 + \nu_1$, $\lambda_1 + \mu_1 + \nu_2$.  Since these line segments need to have non-negative length, we obtain the necessary conditions
$$ \lambda_2 + \mu_1 + \nu_1, \lambda_1 + \mu_2 + \nu_1, \lambda_1 + \mu_1 + \nu_2 \geq 0.$$
These inequalities can be rephrased \eqref{trace-sym} as the statement that the quantities $\lambda_1 - \lambda_2$, $\mu_1 - \mu_2$, $\nu_1 - \nu_2$ form the lengths of a triangle.  The reader may verify from some linear algebra that these conditions are indeed necessary and sufficient for \eqref{lmn}.

\begin{figure}[htbp]
  \begin{center}
    \leavevmode
    \epsfig{file=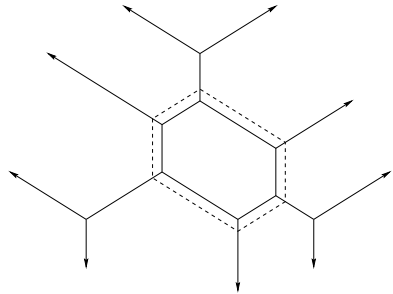,height=1.5in}
    \caption{A hexagon in a honeycomb, with a dotted line indicating a
      place to which one might dilate it.}
    \label{fig:hexbreathe}
  \end{center}
\end{figure}

In the $n > 2$ case things become more complicated, because the boundary values no longer uniquely determine the honeycomb.  In fact, every hexagon which is present in a honeycomb provides a degree of freedom; the hexagon can be ``breathed'' inwards or outwards (see Figure \ref{fig:hexbreathe}).

However, it is still possible to demonstrate that inequalities such
as \eqref{11n} must hold for $n$-honeycombs.  Indeed, one can
simply extend the $\mu_1$ ray northward until it intersects the
$\lambda_1$ ray.  This intersection point must be Northwest of the
intersection of $\nu_n$ and $\lambda_1$, which gives \eqref{11n}.
The reader may be amused by locating proofs of similar inequalities such as $\lambda_1 + \mu_2 + \nu_{n-1} \geq 0$.

\begin{figure}[htbp]
  \begin{center}
   \leavevmode
    \epsfig{file=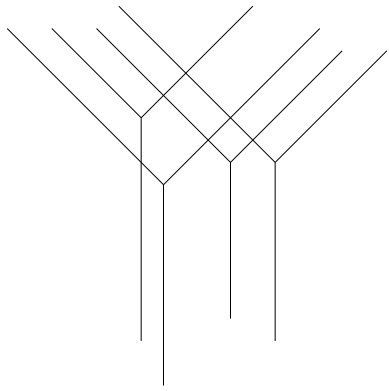,height=2in}
    \caption{An $n$-honeycomb can be obtained by overlaying $n$ 1-honeycombs on top of one another.  Here $\alpha$ and $\beta$ map $1,2,3,4$ to $4, 3, 1, 2$
and $4, 3, 2, 1$ respectively.}
    \label{fig:overlay-1}
  \end{center}
\end{figure}

On the positive side, it is easy to see that \eqref{lmn} will hold if there exist permutations $\alpha, \beta \in S_n$ such that $\lambda_{\alpha(i)} + \mu_{\beta(i)} + \nu_i = 0$ for all $1 \leq i \leq n$.  The honeycomb analogue of this is depicted in Figure \ref{fig:overlay-1}; one can obtain a (rather degenerate) $n$-honeycomb by overlaying $n$ separate 1-honeycombs on top of one another.

\begin{figure}[htbp]
  \begin{center}
    \leavevmode
    \epsfig{file=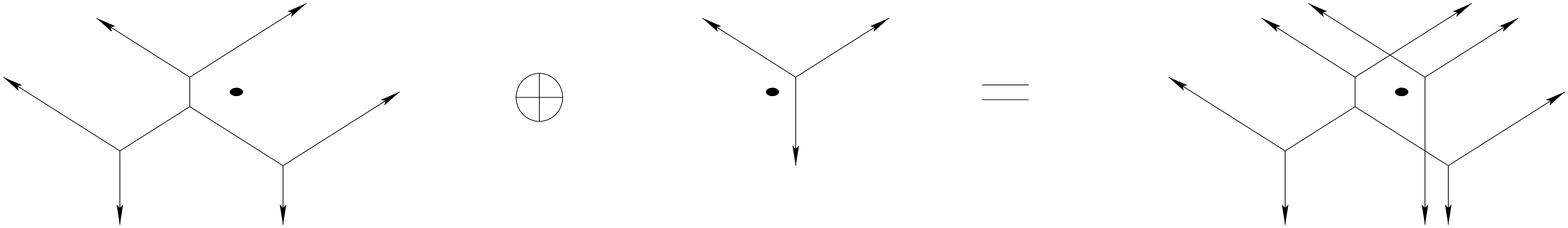,width=6in}
    \caption{Two honeycombs overlaid to produce a third. The origin
      $(0,0,0)$ is marked in each picture with a black dot.}
    \label{fig:overlay}
  \end{center}
\end{figure}

More generally, there is a notion of \emph{overlaying} an $n$-honeycomb $h$ and an $m$-honeycomb $h'$ to form an $n+m$-honeycomb $h \oplus h'$.  To be precise, $h \oplus h'$ is the honeycomb whose associated measure is the sum of the measures associated to $h$ and $h'$.  This operation corresponds to the direct sum operation on Hermitian matrices (which takes an $n \times n$ matrix and an $m \times m$ matrix and forms an $(n+m) \times (n+m)$ block-diagonal matrix), or of the concatenation operation on spectra (which takes a set of $n$ eigenvalues and a set of $m$ eigenvalues and forms the (sorted) set of $n+m$ eigenvalues.)
Intuitively, an overlay can be demonstrated
by drawing two honeycombs on transparencies
and stacking both transparencies on the same projector; 
see figure \ref{fig:overlay}.
We shall have more to say about overlays later in this article.

The relation \eqref{lmn} is clearly $S_3$ symmetric, however the statement that $(\lambda, \mu, \nu)$ admits a honeycomb with these boundary values only appears to be symmetric under cyclic permutations and the inversion $(\lambda, \mu, \nu) \mapsto (-\nu, -\mu, -\lambda)$.  Nevertheless, it is possible to show that honeycombs do indeed obey the full $S_3$ symmetry; an elegant proof of this based on scattering arguments is given in \cite{Wo}.  The same argument also gives the associativity property
\begin{equation}\label{assoc}
(\exists \nu: \lambda \boxplus \mu \sim_c \nu; \nu \boxplus \rho \sim_c \sigma)
\iff (\exists \nu': \mu \boxplus \rho \sim_c \nu'; \lambda \boxplus \nu \sim_c \sigma);
\end{equation}
this property, combined with a ``Pieri rule'' to handle the generating cases when $\lambda$ or $\mu$ is equal to $(\eps, 0, \ldots, 0)$ for some small $\eps$, can be used to give an inductive proof of Theorem \ref{thm:cexist}.

An interesting degenerate case occurs when $\nu$ is kept fixed, while the spacings between eigenvalues of $\lambda$ and $\mu$ are allowed to become very large.  In this case the honeycomb degenerates into a pattern known as a Gelfand-Cetlin pyramid, while Weyl's problem degenerates to Schur's problem of determining which $n$-tuples can be the diagonal entries of a Hermitian matrix with specified eigenvalues.  (In fact, we discovered honeycombs by extrapolating from this degenerate case.)

We encourage the reader to try out the honeycomb Java applet at
\begin{center}
{\tt http://www.math.ucla.edu/\~{}tao/java/honeycombs.html}  
\end{center}

\section{Organizing honeycombs into a polyhedral cone}\label{sec:organize}

The space $\HONEY_n$ of all $n$-honeycombs has been defined as an abstract set, but one can in fact give this space the structure of a polyhedral cone inside some finite dimensional vector space.
 
Call a honeycomb {\dfn nondegenerate} if
\begin{enumerate}
\item all its edges are multiplicity 1, and
\item all its vertices are either right-side-up or upside-down Y's.
\end{enumerate}

\begin{figure}[htbp]
  \begin{center}
    \leavevmode
    \epsfig{file=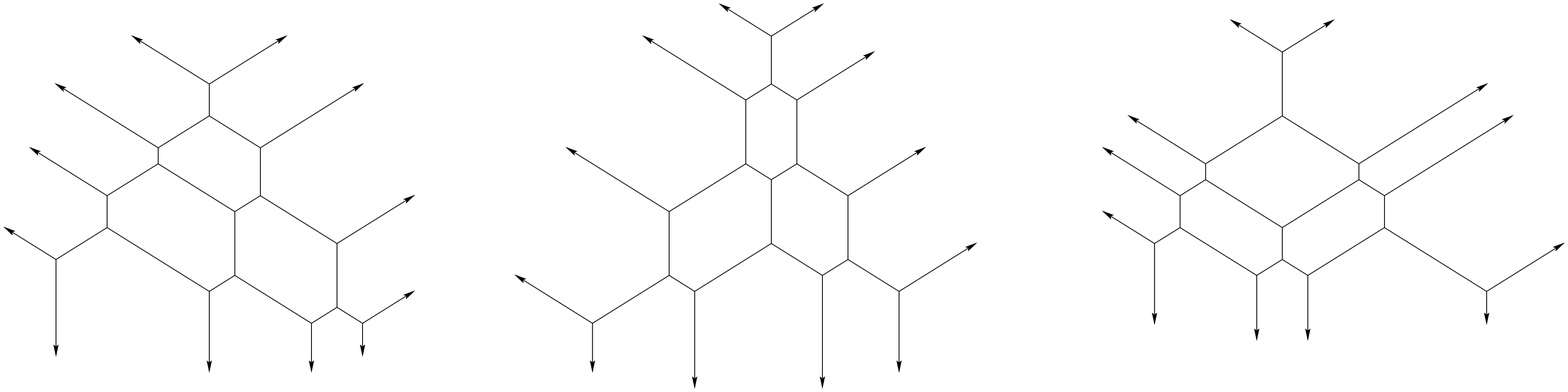,width=6in}
    \caption{Three nondegenerate $3$-honeycombs. Note that there is a
      natural way to correspond the edges in one with the edges in any other.}
    \label{fig:exnondeg}
  \end{center}
\end{figure}

It is straightforward to prove that all nondegenerate $n$-honeycombs have the same topological shape, namely that of \ref{fig:exnondeg}; 
in particular there is a natural way to correspond the edges in one
$n$-honeycomb with those in any other.

This gives us a way of making the space of nondegenerate $n$-honeycombs
into an open polyhedral cone in a real vector space, with defining hyperplanes 
given by rational co-efficients. The
coordinates are given by the constant coordinates of the edges, linear
equalities must be imposed saying that the constant coordinates
of three edges meeting at a vertex add to zero,
and linear inequalities which say that every edge has strictly
positive length.

\begin{Theorem}\cite{KT}\label{thm:picturesdontlie}
  The identification in the above paragraph, between nondegenerate
  $n$-honeycombs and points in a certain rational
  polyhedral open cone, extends to an identification of all $\HONEY_n$
  with the closure of this cone. In particular, nondegenerate $n$-honeycombs 
  form a dense open set in $\HONEY_n$.
\end{Theorem}

This theorem is surprisingly annoying to prove, and takes pp. 1067-1074
of \cite{KT}. Its virtue is in enabling us to use the theory of such cones
to define the special honeycombs that we will deform to, as explained
in the introduction.

For a honeycomb $h$, let $\partial h \in (\reals^n)^3$ denote the list 
$(\lambda_1,\Ldots,\lambda_n,\mu_1,\Ldots,\mu_n,\nu_1,\Ldots,\nu_n)$
of constant coordinates on the boundary edges of $h$, and
$\BDRY_n \leq (\reals^n)^3$ be the image of this map 
$\partial : \HONEY_n \to (\reals^n)^3.$ Then we can think of the
main question as being to list the inequalities determining $\BDRY_n$. 

It is not hard to show directly that this map is proper, so each fiber 
is a compact, convex polyhedron.

As one application of this formalism, we can easily show that for any $\lambda$, $\mu$, $\nu$, the truth or falsity of \eqref{lmn} or \eqref{class-eq} can be determined in polynomial time with respect to the dimension $n$; this fact appears to be previously unknown.  Indeed, the problem is equivalent to determining whether the polytope $\partial^{-1}(\lambda, \mu, \nu)$ is non-empty.  Since the cone $\HONEY_n$ has about $O(n^2)$ faces, this can be achieved in polynomial time by standard algorithms (e.g. the simplex method).  [The authors thank Peter Shor for pointing out this fact].

One can also use this formalism to create a more quantitative version of Theorems \ref{thm:cexist} and \ref{thm:cexist-sym}.  Let $H_\lambda$ be the random variable with the uniform distribution on $\calO_\lambda$ (where ``uniform'' can be defined using induced Lebesgue measure, or the $U(n)$ action).  In other words, $H_\lambda$ is a random matrix with spectrum $\lambda$.  One can then define $P(\lambda \boxplus \mu \sim_c \nu)$ to be the probability density of the spectrum of the sum $H_\lambda + H_\mu$ of two independent random matrices evaluated at $\nu$.  Similarly define $P(\lambda \boxplus \mu \boxplus \nu \sim_c 0)$.

\begin{Theorem}\label{thm:prob}
Up to inessential factors (constants and Vandermonde determinants), $P(\lambda \boxplus \mu \sim_c -\nu)$ and $P(\lambda \boxplus \mu \boxplus \nu \sim_c)$ are equal to the volume of $\partial^{-1}(\lambda, \mu, \nu)$. 
\end{Theorem}

Readers familiar with symplectic geometry will recognize this type of theorem from the theory of moment maps of compact Lie groups such as $U(n)$.  
Indeed, $P(\lambda \boxplus \mu \sim_c -\nu)$ is essentially the volume
of the symplectic reduction of the manifold $\calO_\lambda \times \calO_\mu$
(with the diagonal $U(n)$ action) at the point $-\nu$, and similarly for $P(\lambda \boxplus \mu \boxplus \nu \sim_c 0)$.  We shall have more to say about this later on.

\section{Quantum analogues}\label{quantum-sec}

We now describe the quantum analogue \eqref{quantum-eq} of the classical relation \eqref{class-eq}.  Roughly speaking, \eqref{quantum-eq} is to the representation theory of $U(n)$ as \eqref{class-eq} is to the symplectic geometry of $U(n)$.

Recall that the irreducible unitary representations of $U(1)$ are all one-dimensional.  In fact, for each integer $\lambda$ we can define the irreducible representation $V_{\lambda}$ as a one-dimensional vector space, with the action of $e^{i\theta}$ given by multiplication by $e^{i\lambda\theta}$ on $V_\lambda$. 

More generally, for any weakly decreasing sequence $\lambda = (\lambda_1 \geq \ldots \geq \lambda_n)$ of integers we can define an irreducible unitary representation $V_\lambda$ of $U(n)$ by standard constructions (see e.g. \cite{fulton:young}).  The $n$-tuple $\lambda$ is referred to as the \emph{weight} of $V_\lambda$.  For instance, if $\lambda$ consists of $k$ 1's and $n-k$ 0's then $V_\lambda$ is the space of $k$-forms $\bigwedge^k \complexes^n$ with the standard $U(n)$ action.  More generally, if $\lambda_n \geq 0$ we define $V_\lambda$ to be the highest weight irreducible representation in 
$$ \Tensor_{i=1}^n Sym^{\lambda_i - \lambda_{i+1}} {\bigwedge}^i \complexes^n$$
with the convention $\lambda_{n+1} = 0$, and the $\lambda_n < 0$ representations can be defined dually.

Given two irreducible representations $V_\lambda$, $V_\mu$ of $U(n)$, the tensor product $V_\lambda \otimes V_\mu$ is another representation of $U(n)$.  In the $n=1$ case the tensor product is again an irreducible representation: $V_\lambda \otimes V_\mu \equiv V_{\lambda + \mu}$.  However, in general the tensor product is not irreducible, and splits up as a direct sum of many smaller irreducible representations $V_\nu$.  We can now define the relation \eqref{quantum-eq} as the statement that a copy of $V_\nu$ appears at least once in the tensor product $V_\lambda \otimes V_\mu$.  Note that the quantum relation is only defined for \emph{integral} $\lambda$, $\mu$, $\nu$, whereas the classical relation \eqref{class-eq} is defined for \emph{real} $\lambda$, $\mu$, $\nu$.

There is a close parallel between \eqref{class-eq} and \eqref{quantum-eq}.  For instance, one can obtain the trace identity \eqref{trace} as a necessary condition for \eqref{quantum-eq} by considering the action of the center $U(1)$ of $U(n)$.  One can similarly obtain the necessary condition \eqref{111} by considering the highest weights of the action of a maximal torus $U(1) \times \ldots \times U(1)$ in $U(n)$.  From a more physical viewpoint, one can view the classical problem as a problem of describing how the moments of inertia of bodies in $\complexes^n$ behave under superposition, while the quantum problem is the problem of describing how the spin states of particles in $\complexes^n$ behave under superposition.  (The $n=2$ case is especially interesting to physicists, since $U(2)$ is closely related to $O(3)$.  In this case every representation $V_\nu$ appears at most once in $V_\lambda \otimes V_\mu$ (this corresponds to the fact that 2-honeycombs are determined by their boundary values), and one can parameterize the decomposition explicitly using the Clebsch-Gordan co-efficients.)

This connection between the classical and quantum problems seems to have been noted first in 
\cite{Lidskii} (and in a more general context in \cite{Heckman},
both in 1982), and appears in detail in \cite{Klyachko}; 
the most natural framework for such results is exposed in \cite{LAA}.  Explicitly, the connection is given by

\begin{Theorem}\label{thm:qandc}  Let $\lambda, \mu, \nu$ be weakly decreasing sequences of $n$ integers.
  \begin{enumerate}
  \item (Quantum implies classical.)  If \eqref{quantum-eq} holds, then \eqref{class-eq} holds.
  \item (Classical implies asymptotic quantum.)  Conversely, if \eqref{class-eq} holds, then there exists an integer $N > 0$ such that $N\lambda + N\mu \sim_q N\nu$.  (Here $N\lambda$ is the sequence $(N\lambda_1, \ldots, N\lambda_n)$.)
\end{enumerate}
\end{Theorem}

From this theorem it is natural to phrase

\begin{conjecture}[Saturation conjecture] One can take $N=1$ in the above theorem.  In other words, \eqref{class-eq} and \eqref{quantum-eq} are equivalent for integer $\lambda$, $\mu$, $\nu$.
\end{conjecture}

This conjecture seems to be special to $U(n)$; the na\"\i ve analogue of this conjecture for other Lie groups can be easily shown to be false.  The saturation conjecture is so named because it is equivalent to the set of triples $(\lambda, \mu, \nu)$ obeying \eqref{quantum-eq} being a saturated sub-monoid of $\integers^{3n}$ 

To attack the saturation conjecture using honeycombs we need a quantum analogue of Theorems \ref{thm:cexist}, \ref{thm:cexist-sym}.  We first phrase a symmetric form of \eqref{quantum-eq}.  We say that
\begin{equation}\label{lmn-q}
\lambda \boxplus \mu \boxplus \nu \sim_q 0
\end{equation}
holds if $V_\lambda \otimes V_\mu \otimes V_\nu$ contains a non-trivial $U(n)$-invariant vector.  It is easy to show that \eqref{lmn-q} is equivalent to $\lambda \boxplus \mu \sim_q -\nu$.

A honeycomb is said to be \emph{integral} if its vertices lie on
$\BZ := \B \cap \integers^3$.  Note that the boundary values of an integral honeycomb are necessarily integers.

\begin{Theorem}\label{thm:qexist}
The relationship \eqref{lmn-q} holds if and only if
there exists an integral honeycomb with boundary values $(\lambda,\mu,\nu)$.  As a corollary, the relationship \eqref{quantum-eq} holds if and only if
there exists an integer honeycomb with boundary values $(\lambda,\mu,-\nu)$.  
\end{Theorem}

Note that Theorems \ref{thm:qexist} and \ref{thm:qandc} imply Theorems \ref{thm:cexist} and \ref{thm:cexist-sym}.  

The problem \eqref{quantum-eq} has had a long history, and a solution is given by the famous Littlewood-Richardson rule.  This rule has been formulated in many different ways, most if which involve Young tableaux; a variant due to Berenstein and Zelevinsky can be easily adapted to give Theorem \ref{thm:qexist}.  (Fulton has also shown that this theorem can be proven directly from the Littlewood-Richardson rule.)  Other proofs are known; for instance, one can combine the quantum version of \eqref{assoc} with Pieri's rule for tensoring a $U(n)$ representation with the tautological $C^n$ representation to give an inductive proof of Theorem \ref{thm:qexist}.

A quantum analogue of Theorem \ref{thm:prob} is also known:

\begin{Theorem}\label{q-prob}
The number of times $V_\nu$ appears in the tensor product of $V_\lambda \otimes V_\mu$ is equal to the number of integral honeycombs with boundary values $(\lambda, \mu, -\nu)$.  Equivalently, the dimension of the $U(n)$-invariant subspace of $V_\lambda \otimes V_\mu \otimes V_\nu$ is equal to the number of integral honeycombs with boundary values $(\lambda, \mu, \nu)$.
\end{Theorem}

All the proofs of Theorem \ref{thm:qexist} mentioned above can also be used to prove Theorem \ref{q-prob}.  Theorem \ref{thm:prob} can be viewed as a crude asymptotic version of Theorem \ref{q-prob}.  Variants of this theorem appear in \cite{Johnson,BZ}, and particularly in \cite{GP}, though honeycombs are not explicitly used in these papers.  We remark that the representation theoretic quantities in Theorem \ref{q-prob} can also be calculated by the Steinberg product rule (for instance), though we do not know a proof of this Theorem that goes via this rule.

Readers who are familiar with the representation theory of $SL_2$ (or $SU(2)$) may verify that the honeycomb rule given in Theorem \ref{thm:qexist} corresponds to the usual
triangle inequalities for the weights.  The fact that 2-honeycombs are
uniquely determined by their boundary values corresponds to the fact that
each irreducible representation of $SL_2$ appears exactly once in a tensor
product of irreducibles.

As mentioned in the introduction, the saturation conjecture would give a complete solution to the (now equivalent) problems \eqref{class-eq}, \eqref{quantum-eq}, given by Horn's conjecture.  We shall keep the reader in suspense as to the actual statement of Horn's conjecture for a little while longer, and instead turn to the honeycomb-based proof of saturation.

\section{Proof of the saturation conjecture}\label{-sec}

In light of the theorems of the previous section, the saturation conjecture can be reduced to the following purely honeycomb-theoretic problem:

\begin{Theorem}\label{honey-sat}
Let $h$ be a (real-valued) honeycomb with integer boundary values.  Then there exists an integer honeycomb $h'$ with the same boundary values as $h$.
\end{Theorem}

Or in other words: if $\lambda$, $\mu$, $\nu$ are integers and the polytope $\partial^{-1}(\lambda, \mu, \nu)$ is non-empty, then $\partial^{-1}(\lambda, \mu, \nu)$ must contain at least one integer point.

The most obvious thing to do is to look for a vertex of $\partial^{-1}(\lambda, \mu, \nu)$; however one can give examples of vertices which are non-integral even when $\lambda$, $\mu$, $\nu$ are integers.  Thus we have to be a little more careful as to how to locate our integer honeycomb.

We call a functional 
$f: \HONEY_n \to \reals$ {\dfn superharmonic} if it  
increases when we dilate a hexagon (which one can do to any hexagon in
any nondegenerate honeycomb, as in figure \ref{fig:hexbreathe}.)

Fix a generic superharmonic functional $f$, and define the 
{\dfn largest lift} of a triple
$(\lambda,\mu,\nu)$ as the honeycomb $h$ that maximizes $f(h)$
subject to $\partial h = (\lambda,\mu,\nu)$.
It is straightforward to prove that the largest-lift map
$\BDRY_n \to \HONEY_n$ is uniquely defined (for a given generic $f$), 
continuous, and piecewise-linear.  To show Theorem \ref{honey-sat}, it then suffices to show that the largest map takes integer boundary values to integer honeycombs.

We need some more notation. Say that a honeycomb $h$ has {\dfn only simple degeneracies} if all its edges are multiplicity 1, and its vertices are either Ys (possibly upside down) or crossings of two straight lines. In this case,
define the {\dfn underlying graph of $h$} as the 
graph whose vertices consist of the (possibly upside down) Ys but 
{\em not} the crossings; the crossings we instead interpret as two edges
missing one another. (In the examples in figure \ref{fig:honeyex},
every vertex {\em except} the bottom right vertex of 
the right-hand honeycomb are only simple degeneracies.)

\begin{figure}[htbp]
  \begin{center}
    \leavevmode
    \epsfig{file=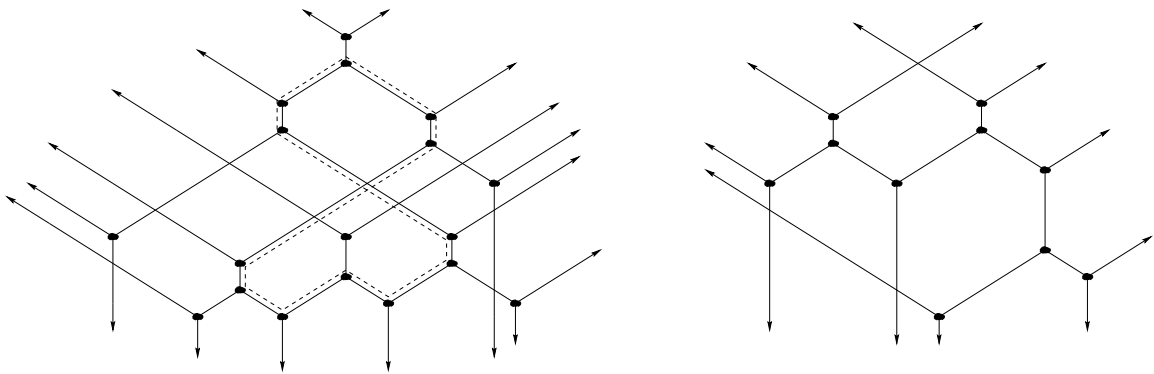,width=6in}
    \caption{Two simply degenerate honeycombs, with black dots on the
      vertices of the their underlying graphs. The left one has a loop
      that can breathe in and out (to, say, the dotted-line position),
      but the right has none.}
    \label{fig:loops}
  \end{center}
\end{figure}

With this in mind, we can talk about {\dfn loops} in 
a simply degenerate honeycomb (meaning, in the underlying graph), 
or call the honeycomb {\dfn acyclic} if there are none. 
For example, in the left honeycomb in figure \ref{fig:loops} there is a loop,
whereas the honeycomb on the right is acyclic.

The importance of loops in simply degenerate honeycombs is that they
can be breathed in and out, as in figure \ref{fig:loops}, 
generalizing the case of dilating a single hexagon. 

Call a largest lift \emph{regular} if the boundary spectra
$\lambda$, $\mu$, $\nu$ each contain no repeated eigenvalues.

The main technical part of \cite{KT} is to prove

\begin{Theorem}\cite{KT}\label{thm:simpdegen}
  Regular largest lifts can only have simple degeneracies.
\end{Theorem}

In particular, regular largest lifts come with underlying graphs.  Roughly speaking, this theorem is proven by showing that every non-simple degeneracy can be ``blown up'' in a way which increases the superharmonic functional.

\begin{Lemma}\label{lem:acyclic}
  The underlying graphs of regular largest lifts are acyclic. 
\end{Lemma}

\begin{proof}[Sketch of proof.]
  If a simply degenerate honeycomb has a loop, we can breathe it in
  and out; one direction will increase the value of any (generic)
  superharmonic functional. A largest lift is by assumption already at
  the maximum value of the functional, so there can be no loops.
\end{proof}

From this lemma one can show that the co-ordinates of regular largest lifts are integral linear combinations of the boundary values.  Those who wish to see the details should go to \cite{KT}, but the
argument is intuitively clear. Given a honeycomb with some edges labeled
by their constant coordinates and some still mysterious, look for vertices
with two known constant coordinates. The remaining one is minus the sum 
of the other two. Label it such and repeat.
The reader is invited to play this game on the 
honeycombs in figure \ref{fig:loops}, and see how in the left-hand
honeycomb one gets stuck exactly because of the loop. 

\begin{figure}[htbp]
  \begin{center}
    \leavevmode
    \epsfig{file=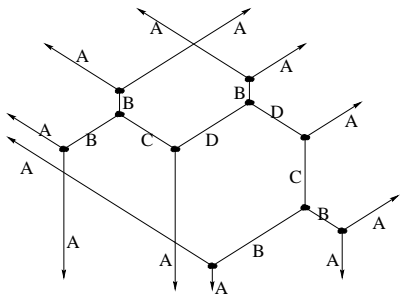,height=1.5in}
    \caption{A honeycomb integrally determined by its boundary,
      in four stages A,B,C,D.}
    \label{fig:intlift}
  \end{center}
\end{figure}

In figure \ref{fig:intlift} we have labeled the boundary edges ``A'',
the edges whose constant coordinates can be determined from those ``B'',
those at the next stage of this recursive algorithm ``C'', and so on.

Since every largest lift can be obtained as a limit of regular largest lifts, we thus have that the co-ordinates of all largest lifts are integral linear combinations of their boundary values.  In particular, if the boundary values are integral, then the largest lifts are also linear.  This proves Theorem \ref{honey-sat}, which gives the saturation conjecture.

\section{Klyachko's result and Horn's conjecture}\label{sec:klyhorn}

In this section we state a version of Klyachko's result (one direction
of which was also proven by Helmke and Rosenthal) and Horn's conjecture.
(We give slightly revisionist versions in order to avoid introducing 
Schubert calculus on Grassmannians, which is one of the equivalent problems 
explained in \cite{FultonBAMS}.)

Horn \cite{Horn} showed that the solution set to \eqref{class-eq} must be given
by \eqref{trace} and a finite number of inequalities of the form
\begin{equation}\label{ijk}
\lambda_{i_1 + r} + \ldots + \lambda_{i_r + 1}
+ \mu_{j_1 + r} + \ldots + \mu_{j_r + 1}
\geq 
\nu_{k_1 + r} + \ldots + \nu_{k_r + 1}
\end{equation}
where $1 \leq r < n$ and $i = (i_1 \geq \ldots \geq i_r)$, $j = (j_1 \geq \ldots \geq j_r)$, and $k = (k_1 \geq \ldots \geq k_r)$ are weakly decreasing sequences of integers between 0 and $n-r$ inclusive.  Let us call triples $(i,j,k)$ of this form \emph{admissible}.

As an example, \eqref{111} is \eqref{ijk} for the admissible triple $((0), (0), (0))$, while Weyl's inequalities correspond to admissible triples of the form $((i), (j), (i+j))$.  The inequality $\lambda_1 + \lambda_2 + \mu_1 + \mu_2 \geq \nu_1 + \nu_2$ corresponds to $((0,0), (0,0), (0,0))$ and so forth.

After many computations of such inequalities, Horn made the remarkable

\begin{conjecture}[Horn's conjecture]  Let $\lambda, \mu, \nu$ be weakly decreasing sequences of reals.  Then \eqref{class-eq} holds if and only if \eqref{trace} holds, and \eqref{ijk} holds whenever $i, j, k$ are admissible and $i \boxplus j \sim_c k$.  
\end{conjecture}

Note that $i, j, k$ have length $r$, which is strictly less than $n$, and that the set of all admissible triples $(i,j,k)$ is finite for each fixed $n$.  Thus, if Horn's conjecture is true, the solution set to \eqref{class-eq} can in principle be determined by a recursive algorithm.  (In fact, Horn's original formulation of this conjecture was phrased in terms of this algorithm.)

Helmke, Rosenthal, and Klyachko showed that Horn's conjecture was true provided that the $\sim_c$ relation on $(i,j,k)$ was replaced by the quantum counterpart $\sim_q$:

\begin{Theorem}\label{cq} Let $\lambda, \mu, \nu$ be weakly decreasing sequences of reals obeying \eqref{trace}.
  \begin{enumerate}
  \item \cite{HR,Klyachko} ] If \eqref{class-eq} holds, then \eqref{ijk} holds whenever $(i,j,k)$ are admissible triples obeying $i \boxplus j \sim_q k$.
  \item \cite{Klyachko} Conversely, if \eqref{ijk} holds whenever $(i, j, k)$ are admissible triples obeying $i \boxplus j \sim_q k$, then \eqref{class-eq} holds.
\end{enumerate}
\end{Theorem}

So in a sense, solvability of the ``classical problem in dimension $n$'' 
(about summing $n\times n$ Hermitian matrices) is determined by the 
solvability of the ``quantum problem in dimension $m<n$'' 
(about tensoring representations of $U(m)$.)  Given the saturation theorem proven in the last section, which says that
each such quantum problem is solvable exactly if the corresponding
classical problem (in the same dimension) is solvable, this gives a
recursive way to answer the problem. 

Theorem \ref{cq} connects the classical and quantum problems in a way markedly different to the standard classical/quantum analogy as codified by Theorem \ref{thm:qandc}.  The proofs of this theorem are highly non-trivial, and first proceed by showing \eqref{quantum-eq} is equivalent to a certain intersection problem in the Schubert calculus of Grassmanians.  We do not discuss this further here, but refer the interested reader to \cite{FultonBAMS}.  More recently, a purely honeycomb-theoretic proof of Theorem \ref{cq} has been obtained, which we discuss briefly in the next section.

\section{Other consequences}\label{sec:othconseq}

We close with mention of a few other applications of honeycombs and
their properties proven above.

Horn's proof that the solution set of \eqref{class-eq} is determined by
\eqref{trace} and a finite number of inequalities of the form \eqref{ijk}
is based on the following stronger fact: if \eqref{ijk} held with equality and $\lambda, \mu, \nu$ are regular, then the associated triple of matrices
$(H_\lambda,H_\mu,H_\nu)$ is necessarily block diagonalizable. Put another way,
the Hermitian triple is the direct sum of two smaller
Hermitian triples.

Given that we have already drawn an analogy between direct sums of matrices and overlaying of honeycombs, there should be a corresponding statement, stating that the faces of $\BDRY_n$ correspond to honeycombs which are overlays. 

Not every direct sum of Hermitian matrices (or an overlay of two honeycombs) corresponds to an inequality \eqref{ijk} in Horn's list.  However, we have

\begin{Theorem}\cite{KTW}
  Let $h$ be a generic point in $\HONEY_n$.  Then $\partial(h)$ is
  on a facet of $\BDRY_n$ if and only if $h$ is an overlay of two smaller
  honeycombs $A$ and $B$, such that at each point of intersection some
  neighborhood of $h$ looks like the left figure 
  in figure \ref{fig:clockwise} (``$A$ turns clockwise to $B$''.)  
  
  \begin{figure}[htbp]
    \begin{center}
      \leavevmode
      \epsfig{file=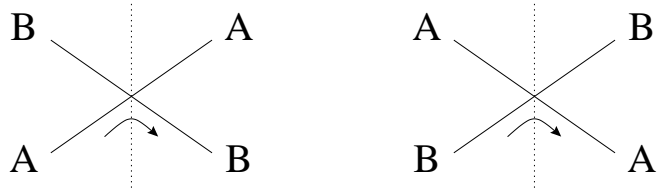,height=1in}
      \caption{In the left figure, $A$ turns clockwise to $B$, whereas in
        the right the reverse is true. Any transverse point of intersection
        of two overlaid honeycombs must look like exactly one of these.}
      \label{fig:clockwise}
    \end{center}
  \end{figure}
\end{Theorem}
  
\begin{figure}[htbp]
  \begin{center}
    \input{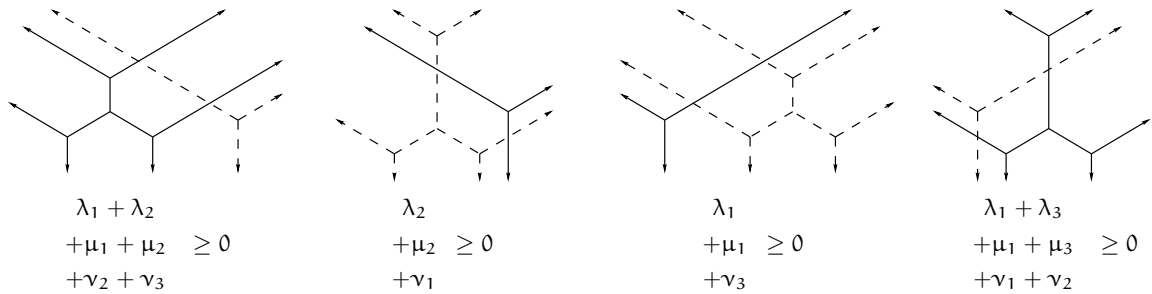}
    \caption{All the consistently clockwise overlays of size $3$,
      up to rotation and deformation, and the associated inequalities
      on $\BDRY_3$.  In each one the solid honeycomb turns clockwise
      to the dashed honeycomb.}
    \label{fig:ovl3}
  \end{center}
\end{figure}

In this case one can read off the inequality on $\BDRY_n$ directly from
$h$: it says that the sum of the boundary coordinates on $A$ are
nonnegative. Some examples are in figure \ref{fig:ovl3}.

\begin{figure}[htbp]
    \begin{center}
      \leavevmode
      \epsfig{file=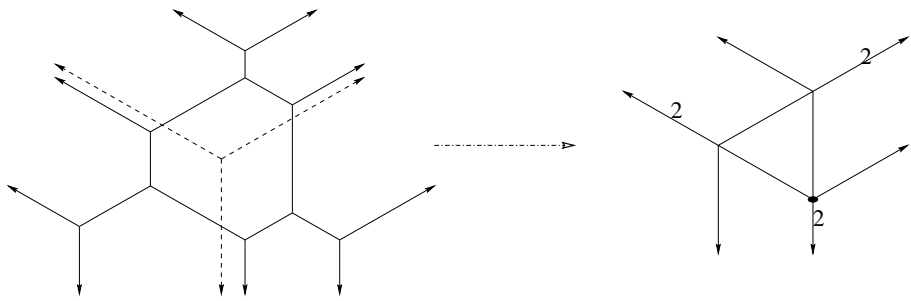,height=1.5in}
      \caption{The honeycomb on the right comes from the solid honeycomb
        on the left, with each edge rescaled to the number of times it
        intersects the dashed honeycomb, and translated to put the
        bottom right vertex at the origin (marked with a dot).}
      \label{fig:shrink4to3}
    \end{center}
  \end{figure}

If $A$ always turns clockwise
to $B$ at the intersections, we can construct a new $m$-honeycomb $A'$
by shifting $A$ so the bottom right vertex is at $(0,0,0)$, replacing
each of $A$'s edges by one whose length is the number of intersections
with $B$, and removing $B$ altogether. The result is in fact a
new honeycomb, integral and of size $m$ 
(an example is in figure \ref{fig:shrink4to3}.)
This construction can be used to give a purely honeycomb-theoretic
proof \cite{KTW} of the results of Klyachko and Helmke-Rosenthal, and gives enough additional insight
to cut down Horn's overcomplete list of inequalities to the 
minimum possible.

\subsection{An open question.} The present proof of Theorem \ref{thm:prob}
is very unsatisfying; it comes as an asymptotic limit of Theorem
\ref{q-prob}, which itself is proved only indirectly.

Consider the horizontal projection of the $2$-sphere of height $1/(2\pi)$
onto the diameter between the poles.  Archimedes' theorem states that
the length of an interval in that diameter equals the area of the
preimage on the sphere. Today we say that the horizontal projection is
{\em measure-preserving}, which at first seems marvelous since the 
interval has only half the dimension of the sphere.  The question: 
is there a corresponding map from $\calM$ to the polytope of honeycombs, 
giving a direct proof of Theorem \ref{thm:prob} (and perhaps of Theorem \ref{q-prob})?

%\goodpagebreak
\bibliographystyle{alpha}

\end{document}